\documentclass[titlepage,11pt]{article}
\usepackage{latexsym}
\usepackage{amsfonts}
\usepackage{amsmath}
\newcommand\blackslug{\hbox{\hskip 1pt \vrule width 4pt height 8pt depth 1.5pt
        \hskip 1pt}}
\newcommand\bbox{\hfill \quad \blackslug \bigbreak}

\def\ll{,\ldots,}
\def\cupcup{\cup\cdots\cup}

\title{Polynomial bounds for chromatic number.\\ I. Excluding a biclique and an induced tree}
\author{Alex Scott\thanks{Research supported by EPSRC grant EP/V007327/1.}\\
Mathematical Institute, University of Oxford, Oxford OX2 6GG, UK
\\
\\
Paul Seymour\thanks{Supported by AFOSR grant
A9550-19-1-0187, and by NSF grant  DMS-1800053.}\\
Princeton University, Princeton, NJ 08544
\\
\\
Sophie Spirkl\thanks{We acknowledge the support of the Natural Sciences and Engineering Research
Council of Canada (NSERC), [funding reference number RGPIN-2020-03912].
Cette recherche a \'et\'e financ\'ee par le Conseil de recherches en sciences
naturelles et en g\'enie du Canada (CRSNG), [num\'ero de r\'ef\'erence
RGPIN-2020-03912].  }\\
University of Waterloo, Waterloo, Ontario N2L3G1, Canada}

\date{}

\newtheorem{thm}{}[section]

\newcommand{\Proof}{\noindent{\bf Proof.}\ \ }

\begin{document}
\maketitle
\begin{abstract}
Let $H$ be a tree. It was proved by R\"odl that graphs that do not contain $H$ as an induced subgraph, and do not
contain the complete bipartite graph $K_{t,t}$ as a subgraph, have bounded chromatic number.
Kierstead and Penrice strengthened this, showing that such graphs have bounded degeneracy.
Here we give a further strengthening, proving that for every tree $H$, the degeneracy is at most polynomial in $t$.  This 
answers a question of Bonamy, Pilipczuk, Rzazewski, Thomass\'{e} and Walczak.
\end{abstract}

\section{Introduction}
The Gy\'arf\'as-Sumner conjecture~\cite{gyarfas, sumner} asserts:
\begin{thm}\label{GSconj}
{\bf Conjecture: } For every forest $H$, there is a function $f$
such that $\chi(G)\le f(\omega(G))$ for every $H$-free graph $G$. 
\end{thm}
(We use $\chi(G)$ and $\omega(G)$ to denote the chromatic
number and the clique number of a graph $G$, and a graph is {\em $H$-free} if it has no induced subgraph isomorphic to $H$.) 
One attractive feature of this conjecture is that it is best possible in a sense: for every graph $H$ that is not a forest, there is no
function $f$ as in \ref{GSconj} (this is easily shown with a random graph).  The conjecture has been proved for some 
special families of trees
(see, for example,~\cite{distantstars, gyarfasprob, gst, kierstead, kierstead2, scott, newbrooms}) but remains open in general.

A class $\mathcal C$ of graphs is {\em $\chi$-bounded} if there is a function $f$ such that $\chi(G)\le f(\omega(G))$ 
for every graph $G$ that is an induced subgraph of a member of $\mathcal C$
 (see~\cite{survey} for a survey).  Thus the
Gy\'arf\'as-Sumner conjecture asserts that, for every forest $H$, the class of all $H$-free graphs is $\chi$-bounded. 
Esperet~\cite{esperet} asked whether, for every $\chi$-bounded class, $f$ can always be chosen to be a polynomial. Neither conjecture has been settled in general.

The complete bipartite graph with parts of cardinality $s,t$ is denoted by $K_{s,t}$.
Let us define $\tau(G)$ to be the largest $t$ such that $G$ contains $K_{t,t}$ as a subgraph (not necessarily induced).
It was proved by R\"odl (mentioned in~\cite{rodl}, and see also~\cite{gst}) that the analogue of the Gy\'arf\'as-Sumner 
conjecture is true if
we replace $\omega(G)$ by $\tau(G)$.  That is:
\begin{thm}\label{hajnalrodl}
For every forest $H$, there is a function $f$
such that $\chi(G)\le f(\tau(G))$ for every $H$-free graph $G$.
\end{thm}
This has the same attractive feature that the result is best possible (in the same sense).

This result was strengthened by Kierstead and Penrice.
Let us say a graph $G$ is {\em $d$-degenerate} (where $d\ge 0$ is an integer) if every nonnull subgraph has a vertex of degree at most $d$;
and the {\em degeneracy $\partial(G)$} of $G$ is the smallest $d$ such that $G$ is $d$-degenerate. Then $\chi(G)\le \partial(G)+1$, 
and so the following result of  Kierstead and Penrice~\cite{kierstead} is a strengthening of \ref{hajnalrodl}:
\begin{thm}\label{kierstead}
For every forest $H$, there is a function $f$
such that $\partial(G)\le f(\tau(G))$ for every $H$-free graph $G$.
\end{thm}

What about the analogue of Esperet's question: do \ref{hajnalrodl} and \ref{kierstead} remain true if we require $f$ to be a polynomial
in $\tau(G)$?
This question was raised by Bonamy, Bousquet, 
Pilipczuk, Rzazewski, Thomass\'{e} and Walczak in~\cite{bonamy}, and they proved it when $H$ is a path, that is:
\begin{thm}\label{bonamy}
For every path $H$, there exists $c>0$ 
such that $\partial(G)\le \tau(G)^c$ for every $H$-free graph $G$.
\end{thm}

In this paper we answer the question completely. Our main result is:
\begin{thm}\label{mainthm}
For every forest $H$, there exists $c>0$
such that $\partial(G)\le \tau(G)^c$ for every $H$-free graph $G$.
\end{thm}

We also look at a related question: what can we say about $\chi(G)$ and $\partial(G)$ if $G$ is $H$-free and does not contain $K_{s,t}$
as a subgraph? More exactly, if $H,s$ are fixed, how do $\chi(G)$ and $\partial(G)$ depend on $t$? We will show
that the dependence is in
fact linear in $t$:
\begin{thm}\label{Kst}
For every forest $H$ and every integer $s>0$, there exists $c>0$ such that for every graph $G$ and every integer $t>0$,
if $G$ is $H$-free and does not contain $K_{s,t}$ as a subgraph, then $\partial(G)\le ct$.
\end{thm}
We also prove a weaker result, that under the same hypotheses, $\chi(G)\le ct$, and for this the bound on $c$ is a small function of $s,H$.

Finally, there is a second pretty theorem in the paper~\cite{bonamy} of Bonamy, Pilipczuk, Rzazewski, Thomass\'{e} and Walczak:
\begin{thm}\label{longhole}
Let $\ell$ be an integer; then there exists $c>0$ such that $\partial(G)\le \tau(G)^c$ for
every graph $G$ with no induced cycle of length at least $\ell$.
\end{thm}
We give a new proof of this, simpler than that in~\cite{bonamy}.

In this paper, all graphs are finite and have no loops or parallel edges. We denote by $|H|$ the number of vertices of a graph $H$.
If $X\subseteq V(G)$, we denote the subgraph of $G$ induced on $X$ by $G[X]$. 
We use 
``$G$-adjacent'' to mean adjacent in $G$, and ``$G$-neighbour'' to mean a neighbour in $G$, and so on.

\section{Producing a path-induced rooted tree.}

We will prove \ref{mainthm} in this section and the next. We need to show that if a graph $G$ has degeneracy at least some very large polynomial in $t$ (independent of $G$), and does not contain $K_{t,t}$
as a subgraph, then it contains any desired tree as an induced subgraph. We will show this in two stages: in this section we will show
that $G$ contains a large (with degrees a somewhat smaller polynomial in $t$) ``path-induced'' tree, and in the next section we will convert 
this to 
the desired induced tree. ``Path-induced'' means that each path of the tree starting at the root is an induced path of $G$; so 
we should be talking about rooted trees. Let us say this carefully.

A {\em rooted tree} $(H,r)$ consists of a tree $H$ and a vertex $r$ of $H$ 
called the {\em root}. A {\rm rooted subtree} of $(H,r)$ means a rooted tree $(J,r)$ where $J$ is a subtree of $H$ and $r\in V(J)$.
The {\em height} of $(H,r)$ is the length (number of edges) of the longest path of $H$ with one end $r$.
If $u,v\in V(H)$ are adjacent and $u$ lies  on the path of $H$ between $v,r$, we say $v$ is a {\em child} of $u$ and $u$ is the {\em parent} of $v$.
The {\em spread} of $H$ is the maximum over all vertices $u\in V(H)$ of the number of children of $u$.
(Thus the spread is usually one less than the maximum degree.)
Let $H$ be a subgraph of $G$ (not necessarily induced), where $(H,r)$ is a rooted tree. We say that
$(H,r)$ is a {\em path-induced rooted subgraph} of $G$ if every path of $H$ with one end $r$ is an induced subgraph of $G$.

Let $\zeta,\eta\ge 1$. The rooted tree $(H,r)$ is {\em $(\zeta,\eta)$-uniform} if
\begin{itemize}
\item every vertex with a child has exactly $\zeta$ children;
\item every vertex with no child is joined to $r$ by a path of $H$ of length exactly $\eta$.
\end{itemize}

We need two lemmas:
\begin{thm}\label{getdisjt}
Let $k,\zeta,\eta\ge 1$ with $\zeta\ge 2$, and let 
$(H_1,r_1)\ll (H_k,r_k)$ be $(k\zeta^{\eta+1},\eta)$-uniform rooted trees, each a subgraph of a graph $G$, such that
$r_i\notin V(H_j)$ for all distinct $i,j\in \{1\ll k\}$. Then for $1\le i\le k$ there is a $(\zeta,\eta)$-uniform rooted 
subtree $(H_i',r_i)$ of $(H_i,r_i)$, such that the trees $H_1'\ll H_k'$ are pairwise vertex-disjoint.
\end{thm}
\Proof
Choose $j\le k$ maximum such that there are  $(\zeta,\eta)$-uniform rooted 
subtrees $(H_i',r_i)$ of $(H_i,r_i)$ for $1\le i\le j$, such that the trees $H_1'\ll H_j'$ are pairwise vertex-disjoint.
Let $X=V(H_1')\cupcup V(H_j')$. Thus $|X|\le j\zeta^{\eta+1}$, since each $H_i'$ has 
$$1+\zeta+\zeta^2+\cdots+\zeta^\eta\le \zeta^{\eta+1}$$
vertices (here we use that $\zeta\ge 2$). Suppose that $j<k$. Then each vertex of $(H_{j+1}, r_{j+1})$ with a child has at least 
$(k-j)\zeta^{\eta+1}\ge \zeta^{\eta+1}\ge \zeta$ children 
not in $X$; and since $r_{j+1}\notin X$, it follows that there
is a $(\zeta,\eta)$-uniform rooted 
subtree $(H_{j+1}',r_{j+1})$ of $(H_{j+1},r_{j+1})$ vertex-disjoint from $X$, contrary to the maximality of $j$. Thus $j=k$, and 
this proves \ref{getdisjt}.~\bbox

Let $t,\eta\ge 1$ and $\zeta\ge 2$ be integers.   Let $(T,r)$ be a $(\zeta,\eta)$-uniform rooted tree, where $T$ is a subgraph of $G$.
A vertex $u\in V(G)\setminus V(T)$ is {\em $t$-bad}
for $(T,r)$ if there is a vertex $w\in V(T)$, with $\zeta$ children in $(T,r)$, such that $u$ is $G$-adjacent to more than $(t-1)\zeta/t$
of these children. 
We will often use the following:
\begin{thm}\label{shrink}
Let $t,\eta\ge 1$ and $\zeta\ge 2$ be integers.   Let $(T,r)$ be a $(t\zeta,\eta)$-uniform rooted tree, where $T$ is a subgraph of $G$;
and let $u\in V(G)\setminus V(T)$. If $u$ is not $t$-bad for $(T,r)$, then there is a $(\zeta,\eta)$-uniform rooted subtree $(S,r)$
of $(T,r)$ such that $u$ has no $G$-neighbour in $V(S)$ except possibly $r$.
\end{thm}
We omit the proof, which is clear.
The second lemma is:

\begin{thm}\label{badverts}
Let $t,\eta\ge 1$ and $\zeta\ge 2$ be integers, where $t$ divides $\zeta$. Let $G$ be a graph that does not contain $K_{t,t}$ as a subgraph, and let $(T,r)$ be a $(\zeta,\eta)$-uniform 
rooted tree, where $T$ is a subgraph of $G$.
Then at most $\zeta^{\eta}(t-1)$ vertices in $V(G)\setminus V(T)$ are $t$-bad for $(T,r)$.
\end{thm}
\Proof
For each $w\in V(T)$ that has $\zeta$ children, let $C_w$ be the set of its children in $(T,r)$. Suppose that
there are $t$ distinct vertices $u_1\ll u_t$ in $V(G)\setminus V(T)$ such that each
is $G$-adjacent to more than $|C_w|(t-1)/t$ vertices in $C_w$. Since $t$ divides $|C_w|$, it follows that each $u_i$ has at most
$|C_w|/t-1$ $G$-non-neighbours in $C_w$, and so at most $t(|C_w|/t-1)$ vertices in $C_w$ are $G$-nonadjacent to one of $u_1\ll u_t$.
Consequently at least $t$ vertices in $C_w$ are $G$-adjacent to all of $u_1\ll u_t$, contradicting that $G$ does not contain $K_{t,t}$ as a subgraph.
Thus there are at most $t-1$ vertices in $V(G)\setminus V(J)$ with more than $|C_w|(t-1)/t$ $G$-neighbours in $C_w$.
So the number of vertices in $V(G)\setminus V(T)$ that are $t$-bad for $(T,r)$ is at most $(t-1)$ times the number of vertices  of $T$
that have children, and so at most $\zeta^{\eta}(t-1)$ (since $\zeta\ge 2$). This proves \ref{badverts}.~\bbox

The main result of this section is the following:
\begin{thm}\label{vertical}
Let $\eta>0$ be an integer and let $c=(\eta+1)!$. Let $\zeta\ge 2$, and let
$(H,r)$ be a rooted tree of height at most $\eta$, and spread at most $\zeta$.
Let $t\ge 1$ be an integer, and suppose that the graph $G$ does not contain $K_{t,t}$ 
as a subgraph, and does not contain a rooted tree isomorphic to $(H,r)$ as a 
path-induced rooted subgraph. Then $\partial(G)\le (\zeta t)^c$.
\end{thm}
\Proof
We may assume that $t\ge 2$.
We proceed by induction on $\eta$. If $\eta=1$, 
it follows that $G$ has maximum degree at most $\zeta-1$, since it does not contain
 $(H,r)$ as a 
path-induced rooted subgraph; and so $\partial(G)\le \zeta-1\le (\zeta t)^c$ as required. So we may assume that $\eta\ge 2$,
and the result holds for all rooted trees with height less than $\eta$. Let $c'=\eta!$
and $\zeta'=t\zeta^{\eta+1}$.
Let us say a {\em limb} is a $(\zeta',\eta-1)$-uniform rooted tree that is a path-induced rooted subgraph of $G$. 
\\
\\
(1) {\em For each vertex $u$, there are at most $\zeta-1$ $G$-neighbours $v$ of $u$ with the property that there is a limb $(J,v)$
of $G$ such that $u\notin V(J)$ and $u$ is not $t$-bad for $(J,v)$.}
\\
\\
Suppose there are $\zeta$ such vertices $v_1\ll v_{\zeta}$, and let the corresponding limbs be $(J_i,v_i)$ for $1\le i\le \zeta$.
By \ref{shrink}, for $1\le i\le \zeta$, there is a $(\zeta^{\eta+1},\eta-1)$-uniform rooted subtree $(J'_i,v_i)$ of $(J_i,v_i)$, 
such that 
$u$ has no neighbour in $V(J'_i)$ except $v_i$.
By \ref{getdisjt}, there is a $(\zeta,\eta-1)$-uniform rooted
subtree $(H_i',r_i)$ of $(J_i',r_i)$ for $1\le i\le \zeta$, such that the trees $H_1'\ll H_k'$ are pairwise vertex-disjoint. But then
adding $u$ to the union of these trees gives a $(\zeta,\eta)$-uniform rooted
tree, and it is path-induced in $G$, and contains a rooted induced subgraph isomorphic to $(H,r)$, a contradiction. This proves (2).

\bigskip 

Let $P$ be the set of vertices $v$ of $G$ such that
there is a limb with root $v$, and let $Q=V(G)\setminus P$. For each $v\in P$, there is at least one limb with root $v$; select one, 
and call it $(J_v,v)$. For each edge $e$ with at least one end in $P$, select one such end, and call it the {\em head} of $e$.
\begin{itemize}
\item Let $A$ be the set of all edges with both ends in $Q$; 
\item Let $B$ be the set of all edges $uv$ of $G$ with head $v$, such that  
$u\notin V(J_v)$, and $u$ is not $t$-bad for $(J_v,v)$;
\item Let $C$ be the set of all edges $uv$ of $G$ with head $v$, such that
$u\notin V(J_v)$, and $u$ is $t$-bad for $(J_v,v)$;
\item Let $D$ be the set of all edges $uv$ of $G$ with head $v$, such that 
$u\in V(J_v)$.
\end{itemize}
Thus every edge of $G$ belongs to exactly one of $A,B,C,D$. 
Since $G[Q]$ does not contain a limb, the inductive hypothesis implies that 
$\partial(G[Q])\le (\zeta' t)^{c'}$. Consequently 
$$|A|\le (\zeta't)^{c'}|Q|\le (\zeta' t)^{c'}|G|.$$
By (1), for each vertex $u\in V(G)$, there are at most $\zeta-1$ edges $uv\in B$ with head $v$; and so 
$$|B|\le (\zeta-1)|G|.$$
For each $v\in P$, there are at most  $\zeta^{\eta-1}(t-1)$ edges $uv\in C$ with head $v$ by \ref{badverts}, and so 
$$|C|\le \zeta^{\eta-1}(t-1)|P|\le \zeta^{\eta-1}(t-1)|G|.$$
For each $v\in P$, since $(J_v,v)$ is path-induced, there are at most $\zeta'$ edges $uv\in D$ with head $v$, and so
$$|D|\le \zeta'|P|\le \zeta'|G|.$$
Summing, we deduce that
$$|E(G)|\le \left((\zeta' t)^{c'}+(\zeta-1)+\zeta^{\eta-1}(t-1)+\zeta'\right)|G|,$$
and so some vertex of $G$ has degree at most $2\left((\zeta' t)^{c'}+(\zeta-1)+\zeta^{\eta-1}(t-1)+\zeta'\right)$.
Since this also holds for every non-null induced subgraph of $G$, we deduce that 
$$\partial(G)\le 2\left((\zeta' t)^{c'}+(\zeta-1)+\zeta^{\eta-1}(t-1)+\zeta'\right).$$
We recall that $\zeta'=t\zeta^{\eta+1}$; and so
$$\partial(G)\le 2\left(\zeta^{c'(\eta+1)} t^{2c'}+(\zeta-1)+\zeta^{\eta-1}(t-1)+\zeta^{\eta+1}t\right).$$
Of the four terms on the right side, the sum of the second and third is at most the fourth, so
$$\partial(G)\le 2\left(\zeta^{c'(\eta+1)} t^{2c'}+2\zeta^{\eta+1}t\right);$$
and since $c'\ge 2$, the second term here is at most the first, so
$$\partial(G)\le 4\zeta^{c'(\eta+1)} t^{2c'}\le \zeta^{c'(\eta+1)} t^{2c'+2}$$
(since we may assume that $t\ge 2$, and so $t^2\ge 4$).
Consequently $\partial(G)\le \zeta^{c'(\eta+1)} t^{2c'+2}.$
But $c=c'(\eta+1)$ and $2c'+2\le c$, and so $\partial(G)\le (\zeta t)^{c}.$ This proves \ref{vertical}.~\bbox

We remark that \ref{vertical} implies \ref{bonamy}, and a strengthening:
\begin{thm}\label{path}
If $H$ is a path, and $t\ge 1$ is an integer, and $G$ is $H$-free and does not contain $K_{t,t}$ as a subgraph, then 
$\partial(G)\le (2t)^{|H|!}$.
\end{thm}
\Proof
Let $\zeta=2$, and $\eta=|E(H)|$. Let $r$ be one end of $H$. Then $G$ does not contain $(H,r)$ as a path-induced rooted subgraph,
and so $\partial(G)\le (2t)^{|H|!}$ by \ref{vertical}. This proves \ref{path}.~\bbox

\section{Growing a tree}

If $(T,r)$ is a rooted tree and $v\in V(T)$, the {\em height} of $v$ in $(T,r)$ is the number of edges in the path between $v,r$; 
and so the 
height of $(T,r)$ is the largest of the heights of its vertices.
Let $(T,r)$ be a rooted tree, and let $(S,r)$ be a rooted subtree. The graph obtained from $T$ by deleting all the edges of $S$
is disconnected, and each of its components contains a unique vertex of $S$; for each $v\in V(S)$, 
let $T_v$ be the component that contains $v\in V(S)$.
We call the rooted tree $(T_v,v)$ the {\em decoration of $S$ at $v$ in $T$}. 

Let $G$ be a graph, let $(S,r)$ be a rooted tree, and let $\zeta\ge 2$ and $\eta\ge 1$. We say that $(S,r)$
is {\em $(\zeta,\eta)$-decorated} in $G$ if $S$ is an induced subgraph of $G$ with height at most $\eta-1$, and there is a rooted tree $(T,r)$ with the following properties:
\begin{itemize}
\item $(S,r)$ is a rooted subtree of $(T,r)$, and $(T,r)$ is a path-induced rooted subgraph of $G$;
\item for each $u\in V(S)$ and $v\in V(T)\setminus V(S)$, if $u,v$ are $G$-adjacent then they are $T$-adjacent;
\item for each $v\in V(S)$, the decoration of $S$ at $v$ in $T$ is $(\zeta,\eta-h)$-uniform, where $h$ is the height of $v$ in $(S,r)$.
\end{itemize}
Thus, informally, $T$ is obtained from $S$ by attaching to $S$ uniform trees rooted at each vertex of $S$. Note that $T$ is only 
required to be path-induced: the various uniform trees that are attached to $S$ might have edges between them.

In view of \ref{vertical}, if we have a graph $G$ with huge degeneracy that does not contain $K_{t,t}$, then it contains
a $(\zeta,\eta)$-uniform rooted tree $(T,r)$ as a path-induced rooted subgraph; and consequently there is a one-vertex
rooted tree $(S,r)$ that is $(\zeta,\eta)$-decorated in $G$. 
The next result shows that if we start with $\zeta$ large enough, then 
by reducing $\zeta$ we can grow $S$ into any larger tree that we wish, and that will prove \ref{mainthm}.

\begin{thm}\label{growtree}
Let $\eta,t\ge 1$ and $\zeta\ge 2$ be integers, let 
$G$ be a graph that does not contain $K_{t,t}$ as a subgraph, and let $(S',r)$ be a $(\zeta',\eta)$-decorated rooted tree in $G$, where
$\zeta'\ge \zeta^{\eta}|S'|t^{\eta+1}$.
Let $p\in V(S')$ with height in $(S',r)$ less than $\eta$. Then there is a $G$-neighbour $q$ of $p$, with $q\in V(G)\setminus V(S')$, 
and with no other $G$-neighbour in $V(S')$,
such that, if $S$ denotes the tree obtained from $S'$ by adding $q$ and the edge $pq$, then $(S,r)$ is a 
$(\zeta,\eta)$-decorated rooted tree in $G$.
\end{thm}
\Proof 
For each $v\in V(S')$, let $h(v)$ denote the height of $v$ in $(S',r)$.
Since  $(S',r)$ is $(\zeta',\eta)$-decorated in $G$, it follows that $S'$ is an induced subgraph of $G$, and there is a rooted tree $(T',r)$ such that
\begin{itemize}
\item $(S',r)$ is a rooted subtree of $(T',r)$, and $(T',r)$ is a path-induced rooted subgraph of $G$;
\item for each $u\in V(S')$ and $v\in V(T')\setminus V(S')$, if $u,v$ are $G$-adjacent then they are $T'$-adjacent;
\item for each $v\in V(S')$, the decoration of $S'$ at $v$ in $T'$ is $(\zeta',\eta-h(v))$-uniform.
\end{itemize}

For each $v\in V(S')$, let $(T_v,v)$ be the decoration of $S'$ at $v$ in $T'$.
Since $T_p$ is $(\zeta',\eta-h(p))$-uniform, and $h(p)<\eta$, it follows that $p$ has $\zeta'$ children in $(T_p,p)$. We need to select 
one of these children, say $q$, to add to $S'$, forming $S$. Any one of them would make a larger induced tree when added to $S'$, since
$(S',r)$ is a
$(\zeta,\eta)$-decorated.
But in order to make the new rooted tree $(S,r)$ $(\zeta,\eta)$-decorated, we will delete from 
$T'$ all vertices of $T'$
that are $G$-adjacent and not $T'$-adjacent to $q$;
and doing so must not destroy too much of $T'$.

For each $v\in V(S')$,  let $(S_v,v)$ be a $(t\zeta,\eta-h(v))$-uniform rooted subtree of $(T_v,v)$. By \ref{badverts}, there are at most
$(t\zeta)^{\eta-h(v)}(t-1)< t^{\eta+1}\zeta^{\eta}$ vertices not in $V(S_v)$ that are $t$-bad for $(S_v,v)$, and so there fewer than 
$\zeta^{\eta}|S'|t^{\eta+1}\le \zeta'$
children of $p$ in $(T_p,p)$ that are $t$-bad for one of the rooted trees $(S_v, v)\;(v\in V(S'))$.
Hence there is at least one child $q$ of
$p$ in $(T_p,p)$ that is $t$-bad for none of the trees $(S_v, v)\;(v\in V(S'))$. Moreover we claim that we can choose $q$ such that $q\notin V(S_p)$.
This is automatic if $(S_p,p)$ has height at least two, since then every child of $p$ in $(S_p,p)$ is bad for $(S_p,p)$, so we may assume
that $(S_p,p)$ has height one, that is, $h(p)=\eta-1$. Consequently no child of $p$ in $(T_p,p)$ is $t$-bad for $(S_p,p)$, and so the 
number that 
are $t$-bad for  one of the rooted trees $(S_v, v)\;(v\in V(S'))$ is at most $\zeta'-t^{\eta+1}\zeta^{\eta}<\zeta'-t\zeta$. This proves
that we can choose $q$ such that $q\notin V(S_p)$.

Let $Q$ be the component containing $q$ of the graph obtained from $T'$ by deleting $V(S)$; thus $(Q,q)$ is 
$(\zeta',\eta-h(p)-1)$-uniform, and so we may choose
a $(\zeta,\eta-h(p)-1)$-uniform rooted subtree $(R_q,q)$ of $(Q,q)$. Note that $q$ has no neighbours in $V(Q)$ except its neighbours in $T'$,
since $(T',r)$ is path-induced.
Since $q$ is not $t$-bad for any of the rooted trees $(S_v,v)\;(v\in V(S'))$, it follows by \ref{shrink} that for each $v$
there is a $(\zeta,\eta-h(v))$-uniform rooted subtree $(R_v,v)$ of $(S_v,v)$ such that $q$ has no $G$-neighbour in $V(R_v)$ except
possibly $v$, and $q$ is $G$-adjacent to $v$ if and only if they are $T'$-adjacent (that is, $v=p$), since $v\in V(S')$ and $(S',r)$ 
is $(\zeta',\eta)$-decorated.
Let $S$ be the tree induced on $V(S')\cup \{q\}$, and 
let $T$ be the union of $T'$, the trees $R_v\;(v\in V(S')\cup \{q\})$ and the edge $pq$.
Then $S$ satisfies the theorem, because the tree $T$ exists. This proves \ref{growtree}.~\bbox

We deduce \ref{mainthm}, which we restate in a strengthened form:
\begin{thm}\label{mainthm2}
Let $\eta,t\ge 1$ and $\zeta\ge 2$. For every rooted tree $(H,r)$ with height at most $\eta$ and spread at most $\zeta$,
let $c=(\eta+3)!|H|$; then
$\partial(G)\le (|H|\zeta t)^c$ for every $H$-free graph $G$ that does not contain $K_{t,t}$ as a subgraph.
\end{thm}
\Proof
It suffices to prove the statement for all trees $H$, and it is helpful to assign a root $s$ to $H$, so $(H,s)$
is a rooted tree. Choose $\eta\ge 1$ and $\zeta\ge 2$ such that $(H,s)$ has height at most $\eta$ and
spread at most $\zeta$. Let $H$ have $k$ vertices. Define $\zeta_k= \zeta$, and for $i=k-1,k-2\ll 1$
let $\zeta_i=i\zeta_{i+1}^{\eta}t^{\eta+1}$.

Let $G$ be an $H$-free graph that does not contain $K_{t,t}$ as a subgraph. 
Suppose that $G$ contains a one-vertex rooted tree that is $(\zeta_1,\eta)$-decorated in $G$.
Choose a maximal rooted subtree $(F,s)$ of $(H,s)$ such that there is a rooted subtree $(S,r)$ of $G$, isomorphic
to $(F,s)$, such that $(S,r)$ is $(\zeta_i,\eta)$-decorated in $G$, where $i=|F|$.
By \ref{growtree}, $i = k$; and so $G$ contains an induced subgraph isomorphic to $H$, a contradiction.

Thus $G$ contains no one-vertex rooted tree that is $(\zeta_1,\eta)$-decorated in $G$. Hence $G$
contains no $(\zeta_1,\eta)$-uniform rooted tree as a path-induced rooted subgraph, and so by \ref{vertical}
(applied with $(H,r)$ replaced by a $(\zeta_1,\eta)$-uniform rooted tree),
$\partial(G)\le (\zeta_1 t)^d$ where $d=(\eta+1)!$.

Now $\zeta_k=\zeta$, and $\zeta_{k-1}=(k-1)\zeta^{\eta}t^{\eta+1}$. For all $i$ with $1\le i\le k-2$, $\zeta_{i+1}\ge it^{\eta+1}$,
and so $\zeta_i=i\zeta_{i+1}^{\eta}t^{\eta+1}\le \zeta_{i+1}^{\eta+1}$. Consequently 
$$\zeta_1\le \zeta_{k-1}^{(k-2)(\eta+1)}\le \left(k\zeta^\eta t^{\eta+1}\right)^{(k-2)(\eta+1)}\le \left(k\zeta t\right)^{(k-2)(\eta+1)^2}.$$
So $\partial(G)\le (k\zeta t)^c$ where
$c=(k-2)(\eta+1)^2(\eta+1)!+(\eta+1)!\le (\eta+3)!k$. This proves \ref{mainthm2}.~\bbox

\section{Excluding $K_{s,t}$}
In this section we prove \ref{Kst}, and before that we prove a weaker statement, with $\partial(G)$ replaced by $\chi(G)$. For the latter we need the following lemma:
\begin{thm}\label{degen}
Let $J$ be a digraph such that every vertex has outdegree at most $k$. Then the undirected graph underlying $J$ has chromatic number at most
$2k+1$.
\end{thm}
\Proof
Let $G$ be the undirected graph underlying $J$. Since every subgraph of $G$ has the property that its edges can be directed so that it has outdegree
at most $k$, it follows that every such subgraph $H$ has at most $k|H|$ edges; and therefore (if it is non-null) has a vertex of degree at most $2k$. Consequently $G$ is $2k$-degenerate, and so is $(2k+1)$-colourable. This proves \ref{degen}.~\bbox

We use \ref{degen} to prove the following (which we include here because the proof gives a relatively small constant $c$, 
although the fact that some $c$ exists follows from \ref{Kst}):
\begin{thm}\label{weakKst}
Let $H$ be a tree and $s\ge 1$ an integer, and let $c=(2s|H|)^{s+|H|}$. Then for every $H$-free graph $G$ and every
integer $t\ge 1$, if $G$ does not contain $K_{s,t}$ as a subgraph then $\chi(G)\le ct$.
\end{thm}
\Proof We will prove this
by induction on $|H|$ (for the same value of $s$). Let $H$ be a tree and $s\ge 0$
an integer, and suppose the theorem holds for all
smaller trees and the same value of $s$. We may assume that $|H|\ge 3$, since the theorem is true for trees with at most two vertices;
let $p\in V(H)$ have degree one, and let $q$ be its $H$-neighbour. Let $H'$ be obtained by deleting $p$ from $H$. Let
$c'=(2s|H'|)^{s+|H'|}$. We observe that
\\
\\
(1) {\em $c\ge \max\left((|H|-2)^{s-1},(s-1)(|H|-2),\left(2(s-2)(|H|-2)+1\right)c'+1\right)$.}
\\

Let $t\ge 1$ be an integer, and let $G$ be an $H$-free graph not containing
$K_{s,t}$ as a subgraph. We will show that $\chi(G)\le ct$. Suppose that this is false, and choose a minimal induced
subgraph $G'$ of $G$ with $\chi(G')>ct$. It follows that every vertex of $G'$ has degree at least $ct$ (since $c$ is an integer).

Let $v\in V(G')$. We say a subset $X\subseteq V(G')\setminus \{v\}$ is a {\em $v$-bag} if there is an isomorphism from $H'$ to
$G[X\cup\{v\}]$
that maps $q$ to $v$. (Thus each $v$-bag has cardinality $|H|-2$.)

Let $v\in V(G')$, and suppose that there are $s-1$ pairwise disjoint $v$-bags, say $X_1\ll X_{s-1}$.
Since $G$ is $H$-free,
every $G$-neighbour $u$ of $v$ either belongs to $X_i$ or has a $G$-neighbour in $X_i$, for $1\le i\le s-1$. In particular,
every $G$-neighbour $u$ of $v$ not in $X_1\cupcup X_{s-1}$ has a $G$-neighbour in each of $X_1\ll X_{s-1}$. But for each choice
of $x_i\in X_i\;(1\le i\le s-1)$ there are at most $t-1$ $G$-neighbours of $v$ $G$-adjacent to each of $x_1\ll x_{s-1}$ (since they are
also all adjacent to $v$, and $G$ has no $K_{s,t}$ subgraph). Consequently there are at most $(t-1)(|H|-2)^{s-1}$ $G$-neighbours of $v$
not in $X_1\cupcup X_{s-1}$; and hence
$$(s-1)(|H|-2)+(t-1)(|H|-2)^{s-1}> ct.$$
Since $ct= c+c(t-1)$, and $(s-1)(|H|-2)< c$, and $(t-1)(|H|-2)^{s-1}\le c(t-1)$, this contradicts (1); so there is no such choice
of $X_1\ll X_{s-1}$.

Choose an integer $r$ maximum such that  there are $r$ pairwise disjoint $v$-bags, say $X_1\ll X_{r}$.
Consequently $r\le s-2$. Let $Y_v=X_1\cupcup X_r$; then from the maximality of $r$, $X\cap Y_v\ne \emptyset$
for every $v$-bag $X$.
Moreover $|Y_v|\le (s-2)(|H|-2)$.

Let $J$ be the digraph with vertex set $V(G')$ in which every vertex in $Y_v$ is $J$-adjacent from $v$, for each $v\in V(G')$.
Thus $J$ has maximum outdegree at most $(s-2)(|H|-2)$, and so by \ref{degen}, the undirected graph $J'$ underlying $J$ has
chromatic number at most $2(s-2)(|H|-2)+1$; and so
$V(G')=V(J')$ can be partitioned into $2(s-2)(|H|-2)+1$ sets each of which is a stable set of $J'$. Let $Z$ be one of these sets.
Then $G[Z]$ is $H'$-free (because otherwise there would be a vertex $v\in Z$, and a subset $X\subseteq Z\setminus \{v\}$,
and an isomorphism from $H'$ to $G[X\cup \{v\}]$ mapping $q$ to $v$, and hence with $X\cap Y_v\ne \emptyset$; but no
vertex of $Y_v$ belongs to $Z$, since $Z$ is stable in $J'$). From the inductive hypothesis,
$\chi(Z)\le c't$, and hence
$$ct<\chi(G)=\chi(G')\le (2(s-2)(|H|-2)+1)c't$$
contrary to (1). This proves \ref{weakKst}.~\bbox

To prove \ref{Kst}, we will need the following strengthening of \ref{kierstead}, also proved in~\cite{kierstead}:
\begin{thm}\label{kierstead2}
For every forest $H$, and every integer $s>0$, there is a tree $S$ such that for every $H$-free graph $G$,
if $G$ contains $S$ as a subgraph,
then $G$ contains $K_{s,s}$ as a subgraph.
\end{thm}
Now we prove \ref{Kst}, which we restate:
\begin{thm}\label{Kst2}
For every forest $H$ and every integer $s>0$, there exists $c>0$ such that for every graph $G$ and every integer $t>0$,
if $G$ is $H$-free and does not contain $K_{s,t}$ as a subgraph, then $\partial(G)< ct$.
\end{thm}
\Proof Let $S$ be as in \ref{kierstead2}, and let $c=|S|^s$; we will show that $c$ satisfies the theorem.
Let $t>0$ be an integer, and let $G$ be an $H$-free graph that does not contain $K_{s,t}$
as a subgraph. Suppose that $\partial(G)\ge  ct$, and choose $G$ minimal with these properties: then every vertex of $G$ has degree
at least $ct$.
\\
\\
(1) {\em Let $R$ be a tree. If every vertex of $G$ has degree at least $t|R|^s$, then $G$ contains a subgraph $T$ isomorphic to $R$,
and $V(T)$ can be ordered as $\{t_1\ll t_n\}$, such that
for $1\le i\le n$, $t_i$ is $G$-adjacent to at most $s-1$ of $t_1\ll t_{i-1}$.}
\\
\\
We prove this by induction on $|R|$. We may assume that $|R|>1$; let $p\in V(R)$ have degree one in $R$, and let $q$ be its
$R$-neighbour. Let $R'$ be obtained from $R$ by deleting $p$. From the inductive hypothesis,
$G$ contains a subgraph $T'$ isomorphic to $R'$,
and its vertex set can be ordered as $\{t_1\ll t_{n-1}\}$, such that
for $1\le i\le n-1$, $t_i$ is $G$-adjacent to at most $s-1$ of $t_1\ll t_{i-1}$. Choose $v\in V(T')$ such that some isomorphism
from $R'$ to $T'$ maps $q$ to $v$. If some $G$-neighbour $u$ of $v$ does not belong to $V(T')$ and has at most $s-1$ $G$-neighbours in
$V(T')$, then we may set $t_n=u$ as required; so we may assume that every $G$-neighbour $u$ of $v$ in $G$ either belongs to $V(T')$
or has at least $s$ $G$-neighbours in
$V(T')$. Let $X\subseteq V(T')$ with $|X|=s$. If there are at least $t$ vertices in $V(G)$ that are $G$-adjacent to every vertex in $X$,
then $G$ contains $K_{s,t}$ as a subgraph, a contradiction. So for each such $X$, there are at most $t-1$
vertices in $V(G)$ that are $G$-adjacent to every vertex in $X$. Since there are most $|R'|^s$ choices of $X$, there are at most
$(t-1)|R'|^s$ vertices in $V(G)\setminus V(T')$ that have at least $s$ $G$-neighbours in $V(T')$. Consequently $v$ has at most $(t-1)|R'|^s$
$G$-neighbours not in $V(T')$. But it has at most $|R'|$ $G$-neighbours in $V(T')$ and so the degree of $v$ in $G$ is at most
$(t-1)|R'|^s +|R'|<t|R|^s$. This proves (1).

\bigskip
Each vertex of $G$ has degree at least $ct=t|S|^s$; let us apply (1) taking $R=S$.  We deduce that
$G$ contains a subgraph $T$ isomorphic to $S$,
and its vertex set can be ordered as $\{t_1\ll t_n\}$, such that
for $1\le i\le n$, $t_i$ is $G$-adjacent in $G$ to at most $s-1$ of $t_1\ll t_{i-1}$. By \ref{kierstead2},
$G[V(T)]$ contains $K_{s,s}$ as a subgraph. Choose $i$ maximum such that $t_i$ belongs to this subgraph;
then $t_i$ is $G$-adjacent to at least $s$ vertices that are earlier in the ordering, a contradiction. This proves \ref{Kst2}.~\bbox

\section{Long holes}

There is another result in the paper by Bonamy et al.~\cite{bonamy}: 
\begin{thm}\label{longhole2}
Let $\ell\ge 2$ be an integer; then there exists $c>0$ that $\partial(G)\le \tau(G)^c$ for
every graph $G$ with no induced cycle of length at least $\ell$.
\end{thm}
In this section we give a simpler proof of this result.

Let $\eta,t\ge 1$ be integers. We say a rooted tree $(H,r)$ is {\em $(t,\eta)$-tapering} if $(H,r)$ has height $\eta$, and every vertex
$v\in V(H)$ of height $i<\eta$ has exactly $t^{\eta-i}$ children. 
For each $v\in V(H)$, let $h(v)$ be its height in $(H,r)$. 

Let $G$ be a graph.
A map $\phi$ from $V(H)$ into $V(G)$ is a {\em $(t,\eta)$-infusion} of $(H,r)$ into $G$ if
\begin{itemize}
\item for all distinct $u,v\in V(H)$, if $u,v\in V(H)$ are $H$-adjacent then $\phi(u),\phi(v)$ are distinct and $G$-adjacent;
\item for each $u\in V(H)$, if $v,w$ are distinct children of $u$ in $(H,r)$, then $\phi(v)\ne \phi(w)$;
\item for every path $P$ of $H$ with one end $r$, the vertices $\phi(v)\;(v\in V(P))$ are all distinct; and 
\item for every path $P$ of $H$ with one end $r$, and for all distinct
$u,v\in V(P)$, $\phi(u),\phi(v)$ are $G$-adjacent if and only if $u,v$ are $H$-adjacent.
\end{itemize}

Let $\phi$ be a $(t,\eta)$-infusion into $G$. We define $V(\phi)=\{\phi(v):v\in V(H)\}$, and we define the {\em root} of $\phi$ to be $\phi(r)$. We say $u\in V(G)$ is {\em $t$-bad}
for $\phi$ if there exists $v\in V(H)$ with $h(v)<\eta$, such that $u$ is distinct from and $G$-adjacent to $\phi(w)$ for more than 
$(t-1)t^{\eta-h(v)-1}$ children $w$ of $v$ in $(H,r)$. Then we have:

\begin{thm}\label{badinfusion}
Let $t,\eta\ge 1$ be integers, let $(H,r)$ be a $(t,\eta)$-tapering rotted tree, let $G$ be a graph not containing $K_{t,t}$ as a subgraph, 
and let $\phi$ be a $(t,\eta)$-infusion of $(H,r)$ into $G$. There are at most $t^{\eta^\eta}$ vertices in $G$ that are $t$-bad
for $\phi$.
\end{thm}
The proof is like that for \ref{badverts}, using that $H$ has 
at most $t^{\eta^\eta-1}$ vertices that have children, and we omit it.

The next result strengthens \ref{longhole}:
\begin{thm}\label{longerhole}
Let $\eta\ge 2$ be an integer, and let $G$ be a graph with no induced cycle of length more than $\eta$. For every integer $t\ge 1$,
if $G$ does not contain $K_{t,t}$ as a subgraph then $\partial(G)\le t^{7\eta^\eta}$.
\end{thm}
\Proof We may assume that $t\ge 2$.
Let $t\ge 1$ be an integer, and let $G$ be a graph with no induced cycle of length more than $\eta$ that does not contain $K_{t,t}$.
Let $(H,r)$ be a $(t,\eta)$-tapering rooted tree (not necessarily contained in $G$). 
\\
\\
(1) {\em If $u\in V(G)$ and $v_i$ is a $G$-neighbour of $u$ for $1\le i\le t^\eta$, all distinct, and for each $i$
there is  a $(t,\eta)$-infusion of $(H,r)$ into $G$ with root $v_i$, such that $u\notin V(\phi_i)$, and $u$ is not $t$-bad
for $\phi_i$, then there is a $(t,\eta)$-infusion of $(H,r)$ into $G$, with root $u$.}
\\
\\
Let $(H',r)$ be a $(t,\eta-1)$-tapering rooted subtree or $(H,r)$.
It follows (analogously to \ref{shrink}) that for $1\le i\le t^\eta$,
there is a $(t,\eta-1)$-infusion $\phi_i'$ of $(H',r)$ into $G$ such that $u$ has no $G$-neighbour in $V(\phi_i')$ except $v_i$. Let us 
number the components of $H\setminus \{r\}$ as $H_1\ll H_{t^\eta}$. Let $\psi(r)=v$, and for $1\le i\le t^\eta$
and each $v\in V(H_i)$, define $\psi(v)=\phi_i'(w)$ where $w$ is the parent of $v$ in $(H,r)$. Then $\psi$ is 
a $(t,\eta)$-infusion of $(H,r)$ into $G$, with root $v$. This proves (1).

\bigskip

In these circumstances we say that $\psi$, constructed as in the proof of (1), is {\em derived from}
the sequence $(\phi_i:\;1\le i\le t^\eta)$.

If $P$ is a path of $H$ with length $\eta$ and one end $r$, and $\phi$ is a $(t,\eta)$-infusion of $(H,r)$ into $G$, then
$\phi$ maps $P$ to an induced path $\phi(P)$ of $G$ with length $\eta$ and with one end the root of $\phi$. We call $\phi(P)$
a {\em column} of $\phi$. We observe that if $\psi$ is derived from $(\phi_i:\;1\le i\le t^\eta)$ as above, then for every column
$Q$ of $\psi$, there is a column $Q'$ of one of $\phi_i(1\le i\le t^\eta)$, say of $\phi'$, such that $Q\setminus \psi(r)$ is a subpath of $Q'$. Let us
call $(\phi',Q')$ a {\em shift} of $(\phi,Q)$.

Let $\mathcal{A}_1$ be the set of all $(t,\eta)$-infusions of $(H,r)$ into $G$. Inductively for $i>1$, let
$\mathcal{A}_i$ be the set of all $(t,\eta)$-infusions $\phi$ such that for some choice of $\phi_i(1\le i\le t^\eta)\in \mathcal{A}_{i-1}$,
$\phi$ is derived from $(\phi_i:\;1\le i\le t^\eta)$. Thus $\mathcal{A}_i\subseteq \mathcal{A}_{i-1}$ for each $i$.
There are two cases: either $\mathcal{A}_i$ is empty for some $i$, or it remains nonempty for all values of $i$.
Suppose first that $\mathcal{A}_i$ is nonempty for all $i$, and let $\mathcal{A}$ be the intersection of all 
the sets $\mathcal{A}_i(i\ge 1)$. Choose $\phi_1\in \mathcal{A}$, and let $Q_1$ be a column of $\phi_1$. 
Since $\phi_1$ is derived from some 
members of $\mathcal{A}$, there exists $\phi_2\in \mathcal{A}$ with root $u_2$, and a column $Q_2$ of $\phi_2$, such that
$(\phi_2,Q_2)$ is a shift of $(\phi_1,Q_1)$. Similarly we can choose an infinite sequence $(\phi_i,Q_i)\;(i = 1,2,3\ldots)$
such that each $\phi_i\in \mathcal{A}$ and each $(\phi,Q_i)$ is a shift of its predecessor. Let $v_i$ be the root of $\phi_i$ for each $i$. Then
$v_i,v_{i+1}\ll v_{i+\eta}$ are the vertices in order of $Q_i$ for each $i$; and so form an induced path of $G$. Since $G$
is finite, there exists $j>0$ such that $v_j$ is adjacent to one of $v_1\ll v_{j-2}$; choose a minimum such value of $j$, and choose $i\le j-2$
maximum such that $v_i,v_j$ are adjacent. Then $\{v_i\ll v_j\}$ induces a cycle of $G$ of length more than $\eta$, a contradiction.

So the second case holds, that is, $\mathcal{A}_i$ is empty for some $i$. Choose $k+1$ minimum such that $\mathcal{A}_{k+1}= \emptyset$.
For $1\le i\le k$ let $X_i$ be the set of all vertices $v$ such that $v$ is the root of a member of $\mathcal{A}_i$ and not the root
of any member of $\mathcal{A}_{i+1}$. Thus the sets $X_1\ll X_k$ are pairwise disjoint. Let $X_0$ be the set of vertices that are
not the root of any member of $\mathcal{A}_1$; so the sets $X_0\ll X_k$ form a partition of $V(G)$.
For each edge $e$ of $G$ with an end in one of $X_1\ll X_k$, choose $i$ maximum such that $e$ has an end in $X_i$, let $v$
be an end of $e$ in $X_i$, and call $v$ the {\em head} of $e$. For each $v\in X_i$, choose $\phi_v\in \mathcal{A}_i$ with root $v$.
(Thus $\phi_v\notin \mathcal{A}_{i+1}$ from the definition of $X_i$.)
Define
\begin{itemize}
\item $A$ is the set of all edges of $G$ with both ends in $X_0$;
\item $B$ is the set of all edges $uv$ with head $v$ such that $u\notin V(\phi_v)$ and $u$ is not bad for $\phi_v$;
\item $C$ is the set of all edges $uv$ with head $v$ such that $u\notin V(\phi_v)$ and $u$ is bad for $\phi_v$;
\item $D$ is the set of all edges $uv$ with head $v$ such that $u\in V(\phi_v)$.
\end{itemize}
Since there is no $(t,\eta)$-infusion of $(H,r)$ into $G[X_0]$, it follows that $G[X_0]$ does not contain
a $(\zeta,\eta)$-uniform tree as a path-induced rooted subgraph, where $\zeta=t^{\eta}$, and so 
$\partial(G[X_0])\le (\zeta t)^{(\eta+1)!}$ from \ref{vertical}. Hence
$$|A|\le (\zeta t)^{(\eta+1)!}|G|.$$
For each $u\in V(G)$, with $u\in X_i$ say, there do not exist $t^\eta$ neighbours $v$ of $u$
such that $uv$ 
has head $v$ and belongs to $B$, since there is no $(t,\eta)$-infusion of $(H,r)$ with root $u$ that is derived from members of
$\mathcal{A}_i$. Hence 
$$|B|\le t^\eta |G|.$$
For each $v\in V(G)$, there are at most $t^{\eta^\eta}$ neighbours $u$ of $v$ such that the edge $uv$ has head $v$ and belongs to $C$,
by \ref{badinfusion}; so 
$$|C|\le t^{\eta^\eta}|G|.$$
Finally, for each $v\in V(G)$, there are at most $t^\eta$ neighbours $u$ of $v$ such that the edge $uv$ has head $v$ and belongs to $D$; so
$$|D|\le  t^\eta|G|.$$
Summing, we obtain
$$|E(G)|\le \left(\left(t^{\eta+1}\right)^{(\eta+1)!}+ t^\eta +t^{\eta^\eta}+ t^\eta\right)|G|\le \left(t^{(\eta+2)!}+t^{\eta^\eta}\right)|G|\le t^{7\eta^\eta}/2.$$
Consequently $\partial(G)\le t^{7\eta^\eta}$. This proves \ref{longerhole}.~\bbox
i
\section*{Acknowledgement}
We would like to express our thanks to Andras Gy\'arf\'as, who clarified the somewhat confusing history of the authorship of \ref{hajnalrodl}
for us.


\begin{thebibliography}{99}
\bibitem{bonamy} 
M. Bonamy, N. Bousquet, M. Pilipczuk, P. Rzazewski, S. Thomass\'{e} and B. Walczak, 
``Degeneracy of $P_t$-free and $C_{\geq t}$-free 
graphs with no large complete bipartite subgraphs'', {\tt arXiv:2012.03686}.
\bibitem{distantstars} 
M. Chudnovsky, A. Scott and P. Seymour, 
``Induced subgraphs of graphs with large chromatic number.  XII. Distant stars'', 
{\em J. Graph Theory} {\bf 92} (2019), 237--254,
{\tt arXiv:1711.08612}.
\bibitem{esperet} L. Esperet,
{\em Graph Colorings, Flows and Perfect Matchings},
Habilitation thesis, Universit\'e Grenoble Alpes (2017), 24.
\bibitem{gyarfas}
A. Gy\'arf\'as, ``On Ramsey covering-numbers'',
in {\em Infinite and Finite Sets, Vol. II} (Colloq., Keszthely, 1973), {\em Coll. Math. Soc. J\'anos Bolyai} {\bf 10}, 801--816.
\bibitem{gyarfasprob}
A. Gy\'arf\'as, 
``Problems from the world surrounding perfect graphs'', 
{\em Proceedings of the International Conference on Combinatorial Analysis and its Applications},  (Pokrzywna, 1985),
{\em Zastos. Mat.}  {\bf 19} (1987), 413--441.
\bibitem{gst}
A. Gy\'arf\'as, E. Szemer\'edi and Zs. Tuza,
``Induced subtrees in graphs of large chromatic number'',
{\em Discrete Math.} {\bf 30} (1980), 235--344.
\bibitem{kierstead}
H. A. Kierstead and S.G. Penrice,
``Radius two trees specify $\chi$-bounded classes'',
{\em J. Graph Theory} {\bf 18} (1994), 119-–129.
\bibitem{rodl}
H. A. Kierstead and V. R\"odl, ``Applications of hypergraph coloring to coloring graphs
not inducing certain trees'', 
{\em Discrete Math.} {\bf 150} (1996), 187--193.
\bibitem{kierstead2} 
H. A. Kierstead and  Y. Zhu, 
``Radius three trees in graphs with large chromatic number'',
{\em SIAM J. Disc. Math.} {\bf 17} (2004), 571--581.
\bibitem{scott}
A. Scott, 
``Induced trees in graphs of large chromatic number'',
{\em J. Graph Theory} {\bf 24} (1997), 297--311.
\bibitem{newbrooms} 
A. Scott and P. Seymour, 
``Induced subgraphs of graphs with large chromatic number.  XIII. New brooms'', 
{\em European J. Combinatorics} {\bf 84} (2020), 103024, {\tt arXiv:1807.03768}.
\bibitem{survey}
A. Scott and P. Seymour,
``A survey of $\chi$-boundedness'', {\em J. Graph Theory} {\bf 95} (2020), 473--504, {\tt arXiv:1812.07500}.
\bibitem{sumner}
D. P. Sumner,
``Subtrees of a graph and chromatic number'', in
{\em The Theory and Applications of Graphs}, (G. Chartrand, ed.),
John Wiley \& Sons, New York (1981), 557--576.

\end{thebibliography}
\end{document}